\documentclass[12pt]{article}
\usepackage{graphicx}
\usepackage{graphics}
\usepackage{amsthm}
\usepackage{amssymb, amsmath}
\title{Integrals in Gradshteyn and Ryzhik:  Hyperbolic and trigonometric functions}
\author{Mark W. Coffey\\
Department of Physics\\
Colorado School of Mines\\
Golden, CO  80401\\
USA\\
mcoffey@mines.edu\\}
\date{February 10, 2018}
\begin{document}
\maketitle
\baselineskip=25 pt
\begin{abstract}
The well known table of Gradshteyn and Ryzhik contains indefinite and definite integrals of
both elementary and special functions.  We give proofs of several entries containing integrands
with some combination of hyperbolic and trigonometric functions.  In fact, we occasionally present
an extension of such entries or else give alternative evaluations.  We develop connections with special cases of special functions including the Hurwitz zeta function. 
Before concluding we mention new integrals coming from the investigation of certain elliptic functions.

\end{abstract}
 
\medskip
\baselineskip=15pt
\centerline{\bf Key words and phrases}
\medskip 

entries of Gradshteyn and Ryzhik, hyperbolic-trigometric integrals, methods of integration, Hurwitz
zeta function, Gamma function, Bessel function

\medskip
\noindent
{\bf 2010 MSC numbers}
\newline{40C10, 33Bxx, 40A10, 11M35}  

\baselineskip=25pt

\pagebreak
\centerline{\bf 1.  Introduction} 
\medskip

After some preliminary considerations,
in this paper we are initially mostly concerned with the entries of sections 4.118--4.119 and 4.121--4.124 of \cite{grad}.
\footnote{We restrict attention to the 1980 4th edition \cite{grad}.
Errata lists in .pdf format may be found at the website http://www.mathtable.com/errata for the
6th, 7th, and 8th editions.  A series of previously published errata appearing in the journal
Math. Comp. may be found at the associated website http://www.mathtable.com/gr/.  We omit this
series of references.}
Specific entries are written in boldface.  Such entries typically have integrands with
hyperbolic and trigonometric functions, sometimes with powers and/or simple rational functions.
\footnote{As discussed in section 5 of this paper, we seem to have found an entry that requires
correction.}
At the end we present examples of integrals with integrands with hyperbolic and trigonometric
functions that have resulted from a research investigation of elliptic functions \cite{armitage,walkerbk}
and their Laurent expansions in terms of Matter and Hurwitz numbers \cite{hurwitz,matter,rieger,mwchnumbers}.

Integrals are of undeniable value to mathematical, scientific, and engineering analyses.
Very select examples are given in \cite{coffeytet,rockmore}.  
In particular, we recall that the denominators of the integrands of the very important Fermi-Dirac
and Bose-Einstein distributions of quantum statistical mechanics may be written in terms of exponential functions and constants or otherwise in terms of hyperbolic functions.
Integrals can provide the solution of differential equations and the explicit formulae and
Voronoi summation \cite{voronoi} of analytic number theory. 
Instead of writing that integrals are irresistible \cite{mollbook}, we may state that they are indispensible.

We have placed some related integrals with hyperbolic-trigonometric denominator integrands in an Appendix.
Some instances of trigonometric-hyperbolic integrals outside of the Gradshteyn and Ryzhik table are 
given in \cite{glasser1,glasser2}.

For some intermediate results it is useful to introduce the analytically continuable Hurwitz zeta
function $\zeta(s,a)=\sum_{n=0}^\infty (n+a)^{-s}$, Re $s>1$, with the obvious functional equation
$\zeta(s,a+1)=\zeta(s,a)-a^{-s}$ \cite{ww}.  The special cases $\zeta(s,1)=\zeta(s)$ and $\zeta(s,1/2)
=(2^s-1)\zeta(s)$ are reductions to the Riemann zeta function $\zeta(s)$.  Explicit evaluations are
$\zeta(0,x)=1/2-x$ and $\zeta'(0,x)=-{1 \over 2}\ln(2\pi)+\ln \Gamma(x)$, with $\Gamma$ the Gamma function.
The Hurwitz zeta function has many known integral representations, and $\zeta(n+1,x)$ for $n >0$
an integer has the reduction to the polygamma function $\psi^{(n)}(x)$ (e.g., \cite{grad}, p. 944):
$$\psi^{(n)}(x)=(-1)^{n+1}n!\zeta(n+1,x).$$

This much sets the stage for a first result which includes many special cases.
{\bf Lemma 1}.  Let Re $p>1$.  Then for $|Re ~b|<|Re ~a|$,
$$\int_0^\infty {{\cosh bx} \over {\sinh ax}}x^{p-1}dx={{\Gamma(p)} \over {(2a)^p}}\left[\zeta\left(p,
{{a-b} \over {2a}}\right)+\zeta\left(p,{{a+b} \over {2a}}\right)\right].$$

{\it Proof}.  Via geometric series expansion of the denominator,
$$\int_0^\infty {{\cosh bx} \over {\sinh ax}}x^{p-1}dx=2\sum_{j=0}^\infty \int_0^\infty x^{p-1} \cosh bx
e^{-(2j+1)ax}dx$$
$$={{\Gamma(p)} \over {(2a)^p}}\sum_{j=0}^\infty \left[{1 \over {(j+(a-b)/2a)^p}}+{1 \over 
{(j+(a+b)/2a)^p}}\right]$$
$$={{\Gamma(p)} \over {(2a)^p}}\left[\zeta\left(p,
{{a-b} \over {2a}}\right)+\zeta\left(p,{{a+b} \over {2a}}\right)\right].$$

{\it Remarks}.  Special cases then include reductions to values of the polygamma functions and/or
the Riemann zeta function.  An example is for $b=0$, when the following representation {\bf 3.523.1} is recovered:
$$\zeta(p)={{(1-2^{-s})^{-1}} \over {2\Gamma(p)}}\int_0^\infty {t^{s-1} \over {\sinh t}}dt, ~~~~
\mbox{Re} ~p>1.$$
Indeed, the proof of Lemma 1 shows {\bf 3.524.5}.  Then with a derivative with respect to $a$ in
{\bf 3.524.3} and the use of analytic continuation, the same result is shown.  In such ways,
{\bf 3.524.1}, {\bf 3.524.3} and {\bf 3.524.5} are equivalent.

Another example is the following, which is equivalent to {\bf 4.111.6}.  We first note the
functional equation $\psi(1-z)-\psi(z)=\pi \cot \pi z$, giving $\psi'(1-z)+\psi'(z)=\pi^2 \csc^2 (\pi z)$,
where $\psi=\Gamma'/\Gamma$ is the digamma function and $\psi'$ is the trigamma function.
{\newline \bf Corollary 1}.
$$\int_0^\infty {{\cosh bx} \over {\sinh ax}}xdx={\pi^2 \over {4a^2}}{1 \over {\sin^2\left({\pi \over 2}
\left(1+{b \over a}\right)\right)}}={\pi^2 \over {4a^2}}{1 \over {\cos^2(\pi b/(2a))}}.$$

{\it Proof}.  From Lemma 1 we have
$$\int_0^\infty {{\cosh bx} \over {\sinh ax}}xdx={1 \over {4a^2}}
\left[\zeta\left(2,{{a-b} \over {2a}}\right)+\zeta\left(2,{{a+b} \over {2a}}\right)\right]$$
$$={1 \over {4a^2}}\left[\psi'\left({1 \over 2}-{b \over {2a}}\right)+\psi'\left({1 \over 2}+{b \over {2a}}\right)\right]$$
$$={\pi^2 \over {4a^2}}{1 \over {\sin^2\left({\pi \over 2}\left(1+{b \over a}\right)\right)}}.$$

There are multiple connections between Lemma 1 and other entries of section 4.111, but we omit
these at this time.  In fact, the result of Lemma 1, and its various derivatives, provide a
starting point for any number of other integrals, by way of power series summation over $p$, with
various coefficients.  This is illustrated specifically in section 3, using further properties of
the Hurwitz zeta and Gamma functions.  We do not pursue it here,
but Lemma 1 and its analogous results then provide a method, using appropriate limits as necessary,
to prove the entries of section 3.525.


\medskip
\centerline{\bf 2.  Concerning sections 4.118, 4.119, 4.121, and 4.122}
\medskip

For these sections, a very useful method is to employ geometric series expansions of the
csch $x$ and sech $x$ factors of the integrands.
A useful preliminary summation, which is an example to be used with little further comment, is the 
following.

{\bf Lemma 2}.
$$\ln \cosh z=\sum_{k=0}^\infty \ln\left(1+{{4z^2} \over {(2k+1)^2\pi^2}}\right)
=\sum_{k=0}^\infty \ln\left(1+{{2iz} \over {(2k+1)\pi}}\right)\left(1-{{2iz} \over {(2k+1)\pi}}\right),$$
and
$$\ln \sinh z-\ln z = \sum_{k=1}^\infty \ln\left(1+{z^2 \over {k^2\pi^2}}\right)
= \sum_{k=1}^\infty \ln\left(1+{{iz} \over {k\pi}}\right)\left(1-{{iz} \over {k\pi}}\right).$$

{\it Proof}.  This follows, for instance from (\cite{grad} 1.521.1 and 1.521.2) for $\ln \cos x$ and
$\ln \sin x-\ln x$.  The origin of these expressions is the Weierstrass factorizations of $\sin z$
and $\cos z$ over their zeros (\cite{grad} 1.431.1 and 1.4313).

{\bf Lemma 3}. (a) For $|\mbox{Im}~b|\leq$ Re $a$ and Re $a>0$,
$$\int_0^\infty e^{-ax}{{\sin bx} \over x}dx=\tan^{-1}\left({b \over a}\right).$$
(b) For Re $a>0$,
$$\int_0^\infty e^{-ax}{{\sin^2 bx} \over x}dx={1 \over 4}\ln \left(1+{{4b^2} \over a^2}\right).$$

{\it Proof}.  (a) The Laplace transform
$$\int_0^\infty e^{-ax}\sin b x ~dx={b \over {a^2+b^2}}$$
is integrated with respect to $a$. (b) The transform
$$\int_0^\infty e^{-ax}\sin^2 b x ~dx={{2b^2} \over {a^3+4ab^2}}$$
is also integrated with respect to $a$.

The next preliminary result affords a proof of {\bf 4.111.7}, and also gives a summation identity for 
the arctangent function.
{\newline{\bf Lemma 4}}.
$$\int_0^\infty {{\sin a x} \over {\cosh \beta x}}{{dx} \over x}=2\tan^{-1}\left(e^{\pi a/2\beta}\right)
-{\pi \over 2}.$$

{\it Proof}.
$$\int_0^\infty {{\sin a x} \over {\cosh \beta x}}{{dx} \over x}=2\sum_{j=0}^\infty (-1)^j 
\int_0^\infty {{\sin a x} \over x}e^{-(2j+1)\beta x}dx$$
$$=2\sum_{j=0}^\infty (-1)^j \tan^{-1}\left({a \over {\beta(2j+1)}}\right)$$
$$=2\sum_{j=0}^\infty (-1)^j {1 \over {\beta(2j+1)}}\int_0^a {{dt} \over {1+[t/(\beta(2j+1))]^2}}$$
$$=2\int_0^a \mbox{sech}\left({{\pi t} \over {2\beta}}\right)dt$$
$$=2\tan^{-1}\left[\tanh\left({{a\pi} \over {4\beta}}\right)\right]$$
$$=2\tan^{-1}\left(e^{\pi a/2\beta}\right)-{\pi \over 2}.$$

{\it Remark}.  We have the alternative expression for the arctangent function
$$\tan^{-1}z={1 \over {2i}}\ln\left({{1+iz} \over {1-iz}}\right). \eqno(2.1)$$

{\bf 4.118}
$$\int_0^\infty {{x\sin ax} \over {\cosh^2 x}}dx=4\int_0^\infty {{x\sin ax} \over {(e^x+e^{-x})^2}}dx$$
$$=4\int_0^\infty e^{-2x} {{x\sin ax} \over {(1+e^{-2x})^2}}dx$$
$$=-4\sum_{m=1}^\infty (-1)^m m \int_0^\infty xe^{-2mx}\sin ax dx$$
$$=-16a\sum_{m=1}^\infty (-1)^m {m^2 \over {(4m^2+a^2)^2}}$$
$$=-16a\sum_{m=1}^\infty (-1)^m\left[{1 \over {4m^2+a^2}}-{a^2 \over {(4m^2+a^2)^2}}\right]$$
$$={\pi \over 4}\left[-2+a\pi \coth\left({{\pi a} \over 2}\right)\right]\mbox{csch}\left({{a \pi} \over 2}
\right)$$
$$=-{d \over {da}}{{\pi a} \over {2\sinh(\pi a/2)}}.$$

{\bf 4.119}
$$\int_0^\infty {{1-\cos px} \over {\sinh qx}}{{dx} \over x}
=2\sum_{j=0}^\infty \int_0^\infty e^{-(2j+1)q x}(1-\cos px){{dx} \over x}$$
$$=4\sum_{j=0}^\infty \int_0^\infty e^{-(2j+1)q x}\sin^2\left({{px} \over 2}\right){{dx} \over x}$$
$$=\sum_{j=0}^\infty \ln\left(1+{p^2 \over q^2}{1 \over {(2j+1)^2}}\right)$$
$$=\ln \cosh{{p \pi} \over {2q}},$$
wherein Lemmas 2 and 3 were applied.

{\bf 4.121.1} This is the difference of two instances of Lemma 4,
$$\int_0^\infty {{\sin ax-\sin bx} \over {\cosh \beta x}}{{dx} \over x}
=2\sum_{j=0}^\infty (-1)^j \int_0^\infty e^{-(2j+1)\beta x}(\sin ax-\sin bx){{dx} \over x}$$
$$=2\sum_{j=0}^\infty (-1)^j \left[\tan^{-1}\left({a \over {\beta(2j+1)}}\right)-\tan^{-1}\left({b \over {\beta(2j+1)}}\right)\right]$$
$$=2\tan^{-1}{{\exp {{a\pi} \over {2\beta}}-\exp {{b\pi} \over {2\beta}}} \over 
{1+\exp{{(a+b)\pi} \over {2\beta}} }}.$$
The last expression follows from the addition formula for the tangent function.

{\bf 4.121.2}
$$\int_0^\infty {{\cos ax-\cos bx} \over {\sinh \beta x}}{{dx} \over x}
=2\sum_{j=0}^\infty \int_0^\infty e^{-(2j+1)\beta x}(\cos ax-\cos bx){{dx} \over x}$$
$$=\sum_{j=0}^\infty \left[\ln\left(1+{b^2 \over \beta^2}{1 \over {(2j+1)^2}}\right)
-\ln\left(1+{a^2 \over \beta^2}{1 \over {(2j+1)^2}}\right)\right]$$
$$=\ln \cosh{{b \pi} \over {2\beta}}-\ln \cosh{{a \pi} \over {2\beta}}.$$

{\it Remarks}.  The results for 4.121.1 and 4.121.2 show the evident antisymmetry of the 
integrals in $a$ and $b$.  It is obvious that two instances of 4.119 may be subtracted
from one another in order to obtain 4.121.2.

{\bf 4.122.1}
This result follows from the trigonometric identity
$$\cos \beta x\sin \gamma x={1 \over 2}[\sin(\gamma+\beta)-\sin(\beta-\gamma)]$$
and the immediate application of 4.121.1:
$$\int_0^\infty {{\cos \beta x\sin \gamma x} \over {\cosh \delta x}}{{dx} \over x}
=\tan^{-1}{{\sinh \gamma \pi/2\delta} \over {\cosh \beta \pi/2\delta}}.$$

{\bf 4.122.2}  By using Lemma 3(b),
$$\int_0^\infty \sin^2 ax {{\cosh \beta x} \over {\sinh x}}{{dx} \over x}$$
$$=\sum_{j=0}^\infty \int_0^\infty \sin^2 ax{{\cosh \beta x} \over x}e^{-(2j+1)x}dx$$
$$={1 \over 4}\sum_{j=0}^\infty \left[\ln\left(1+{{4a^2} \over {(-\beta+(2j+1))^2}}\right)
+\ln\left(1+{{4a^2} \over {(\beta+(2j+1))^2}}\right)\right]$$
$$={1 \over 4}\ln {{\cosh 2a\pi+\cos \beta \pi} \over {1+\cos \beta \pi}}.$$

{\it Remark}.  One could also proceed via the following steps, but the above approach seems
preferable.
$$\int_0^\infty \sin^2 ax {{\cosh \beta x} \over {\sinh x}}{{dx} \over x}$$
$$=\sum_{j=0}^\infty \int_0^\infty (1-\cos2ax){{\cosh \beta x} \over x}e^{-(2j+1)x}dx$$
$$={1 \over 4}\sum_{j=0}^\infty [-2\ln(\beta-1-2j)-2\ln(\beta+1+2j)
+\ln(2ia+(\beta-1-2j))+\ln(2ia-(\beta-1-2j))$$
$$+\ln(2ia+(\beta+1+2j))+\ln(2ia-(\beta+1+2j))].$$

\medskip
\centerline{\bf 3. Elaboration:  Alternative methods for {\bf 4.122.2}}
\medskip

There are of course many relations between the various entries of \cite{grad}.  In this case,
{\bf 3.981.5}, for $|Re ~\beta| <|Re ~\gamma|$,
$$\int_0^\infty \cos ax {{\sinh \beta x} \over {\sinh \gamma x}}dx={\pi \over {2\gamma}}
{{\sin \pi\beta/\gamma} \over {\cosh a\pi/\gamma+\cos \beta \pi/\gamma}}.$$
Since $\int \sinh \beta x ~d\beta={1 \over x}\cosh \beta x$, with the relabelling of variables 
and the use of 4.119 we may recover 4.122.2.

The following is a more involved method based upon the use of Lemma 1 and the next Lemma.
\newline{\bf Lemma 5}.  Letting Re $a>0$,
$$\sum_{n=1}^\infty {z^n \over n}\zeta(2n,a)=-2\ln \Gamma(a)+\ln\Gamma(a-\sqrt{z})+\ln\Gamma(a+\sqrt{z}).
$$

{\it Proof}.  This result may be proved by inserting an integral representation for the Hurwitz
zeta function, or else by first rewriting $\zeta(2n,a)$ in terms of polygamma functions, and 
then inserting an integral representation for the latter.

{\it Other proof for {\bf 4.122.2}}.  By making use of $\sin^2 ax=(1-\cos 2ax)/2$, we have
$$\int_0^\infty \sin^2 ax {{\cosh \beta x} \over {\sinh x}}{{dx} \over x}
={1 \over 2}\int_0^\infty (1-\cos 2ax){{\cosh \beta x} \over {\sinh x}}{{dx} \over x}$$
$$=-{1 \over 2}\sum_{n=1}^\infty {{(-1)^n} \over {(2n)!}}(2a)^{2n} \int_0^\infty x^{2n-1}
{{\cosh \beta x} \over {\sinh x}}dx$$
$$=-{1 \over 4}\sum_{n=1}^\infty {{(-1)^n} \over n}a^{2n}\left[\zeta\left(2n,{{1-\beta} \over 2}\right)
+ \zeta\left(2n,{{1+\beta} \over 2}\right)\right],$$
wherein Lemma 1 has been applied.  Next, Lemma 5 is used twice to give
$$\int_0^\infty \sin^2 ax {{\cosh \beta x} \over {\sinh x}}{{dx} \over x}
={1 \over 2}\ln\left[\Gamma\left({{1-\beta} \over 2}\right)\Gamma\left({{1+\beta} \over 2}\right)\right]$$
$$-{1 \over 4}\ln\left[\Gamma\left({{1-2ia-\beta} \over 2}\right)\Gamma\left({{1+2ia-\beta} \over 2}\right)
\Gamma\left({{1-2ia+\beta} \over 2}\right)\Gamma\left({{1+2ia+\beta} \over 2}\right)\right].$$
Next, we use the property $\Gamma(z)\Gamma(1-z)=\pi/\sin \pi z$ three times, and it is seen that the
terms ${1 \over 2}\ln \pi-{1 \over 4}\ln \pi^2=0$.  Then with $\sin(\pi/2+y)=\cos y$, we arrive at
$$\int_0^\infty \sin^2 ax {{\cosh \beta x} \over {\sinh x}}{{dx} \over x}
=-{1 \over 2}\ln\cos {{\pi \beta} \over 2}+{1 \over 4}\left[\cos\left({{(2ia+\beta)\pi} \over 2}\right)
\cos\left({{(-2ia+\beta)\pi} \over 2}\right)\right].$$
A half-angle formula may be applied to the first term on the right side, and the identity
$$\cos v\cos u={1 \over 2}[\cos(u+v)+\cos(v-u)],$$
to the second term,
$$\cos\left({{(2ia+\beta)\pi} \over 2}\right) \cos\left({{(-2ia+\beta)\pi} \over 2}\right)
={1 \over 2}[\cos(\pi \beta)+\cosh(2a\pi)],$$
since $\cos(ix)=\cosh x$.  Thus we recover the result of 4.122.2, having in the process given
many equivalent expressions.

\medskip
\centerline{\bf 4.  Concerning section 4.123}
\medskip

We begin with a couple of Lemmas that are useful for the evaluation of several entries.

{\bf Lemma 6}. (a)
$${1 \over {\cosh ax+\cos bx}}=-{2 \over {\sin bx}}\sum_{n=1}^\infty (-1)^n e^{-anx}\sin nbx,$$
(b)
$${1 \over {\cosh x-\cos x}}={2 \over {\sin x}}\sum_{n=1}^\infty e^{-nx}\sin nx,$$
and (c)
$$\int_0^\infty {{f(x)dx} \over {\cosh x-\cos x}}=2 \sum_{n=1}^\infty {1 \over n}\int_0^\infty {{f(x/n)} \over {\sin(x/n)}}e^{-x}\sin x dx.$$

{\it Proof}.  We largely omit the proof, but note that for (a) $\sin z={1 \over {2i}}(e^{iz}-e^{-iz})$ so
that geometric series may be applied.  (b) follows from (a) upon replacing $ax$ by $ax-i\pi$.
Then (c) follows for integrable functions $f$. \qed

We note the classical partial fraction expansions of trigonometric functions.  An example is given
next, and such relations will be used without further mention.
\newline{{\bf Lemma 7}.}  
$$\csc \pi x={1 \over {\pi x}}+{{2x} \over \pi}\sum_{k=1}^\infty {{(-1)^k} \over {x^2-k^2}}.$$

\medskip
{\bf 4.123.1}
\medskip

Let Re $a>0$.  Then
$$\int_0^\infty {{\sin x} \over {(\cosh ax+\cos x)}}{{x dx} \over {(x^2-\pi^2)}}$$
$$=-2\sum_{n=1}^\infty \int_0^\infty (-1)^n e^{-anx} \sin nx {x \over {(x^2-\pi^2)}}dx$$
$$=-2\sum_{n=1}^\infty \int_0^\infty (-1)^n e^{-ax} \sin x {x \over {(x^2-n^2\pi^2)}}dx$$
$$=-\int_0^\infty e^{-ax} \sin x\left(\csc x-{1 \over x}\right)dx$$
$$-\int_0^\infty e^{-ax}\left(1-{{\sin x} \over x}\right)dx$$
$$=\tan^{-1}\left({1 \over a}\right)-{1 \over a}.$$

By the same use of Lemmas 3--5 we obtain further examples for Re $a>0$:
$$\int_0^\infty {{\sin 2x} \over {(\cosh ax+\cos 2x)}}{{x dx} \over {(x^2-\pi^2)}}$$
$$=\tan^{-1}{2 \over a}-{{2a} \over {1+a^2}},$$
$$\int_0^\infty {{\sin 3x} \over {(\cosh ax+\cos 3x)}}{{x dx} \over {(x^2-\pi^2)}}$$
$$=\tan^{-1}{3 \over a}-{{4+3a^2} \over {a(4+a^2}},$$
and 
$$\int_0^\infty {{\sin 4x} \over {(\cosh ax+\cos 4x)}}{{x dx} \over {(x^2-\pi^2)}}$$
$$=\tan^{-1}{4 \over a}-{{4a(5+a^2)} \over {9+10a^2+a^4}}.$$

\medskip
{\bf 4.123.2}
\medskip

Let Re $a>0$.  Then
$$\int_0^\infty {{\sin x} \over {(\cosh ax-\cos x)}}{{x dx} \over {(x^2-\pi^2)}}$$
$$=2\sum_{n=1}^\infty \int_0^\infty e^{-ax} \sin x {x \over {(x^2-n^2\pi^2)}}dx$$
$$=\int_0^\infty e^{-ax} \sin x\left(\cot x-{1 \over x}\right)dx$$
$$={a \over {1+a^2}}-\tan^{-1}\left({1 \over a}\right).$$

\medskip
{\bf 4.123.3}
\medskip

Let Re $a>0$.  Then
$$\int_0^\infty {{\sin 2x} \over {(\cosh 2ax-\cos 2x)}}{{x dx} \over {(x^2-\pi^2)}}$$
$$=2\sum_{n=1}^\infty \int_0^\infty e^{-2ax} \sin 2x {x \over {(x^2-n^2\pi^2)}}dx$$
$$=2\int_0^\infty e^{-2ax} \sin 2x\left(\cot x-{1 \over x}\right)dx$$
$$={1 \over {2a}}{{(1+2a^2)} \over {(1+a^2)}}-\tan^{-1}\left({1 \over a}\right).$$

With $\cosh^2 ax-\cos^2 x={1 \over 2}(\cosh 2ax-\cos 2x)$, the addition of 4.123.1 and 4.123.2 gives
{\bf 4.123.4}:  
$$\int_0^\infty {{\cosh ax\sin x} \over {(\cosh 2ax-\cos 2x)}}{{x dx} \over {(x^2-\pi^2)}}$$
$$=-{1 \over {2a(1+a^2)}}.$$
The subtraction of 4.123.1 from 4.123.2 gives 4.123.3.  Thus the pairs (4.123.1,4.123.2) and 
(4.123.3,4.123.4) are equivalent.

\medskip
{\bf 4.123.5}
\medskip
This entry states that for $0<Re ~\beta<1$, Re $\gamma>0$, and $a>0$,
$$\int_0^\infty {{\cos a x} \over {(\cosh \pi x+\cos \pi \beta)}}{{dx} \over {x^2+\gamma^2}}
={{\pi e^{-a\gamma}} \over {2\gamma(\cos \gamma \pi+\cos\beta \pi)}}$$
$$+{1 \over {\sinh \beta \pi}}\sum_{k=0}^\infty\left[{e^{-(2k+1-\beta)a} \over {\gamma^2-(2k+1-\beta)^2}}
-{e^{-(2k+1+\beta)a} \over {\gamma^2-(2k+1+\beta)^2}}\right].$$

Letting
$$I_5 \equiv \int_0^\infty {{\cos a x} \over {(\cosh \pi x+\cos \pi \beta)}}{{dx} \over {x^2+\gamma^2}},$$
we have the differential equation
$${{\partial^2 I_5} \over {\partial a^2}}-\gamma^2 I_5=-\int_0^\infty {{\cos a x} \over {\cosh \pi x+\cos \pi \beta}}dx$$
$$={{2\pi} \over {\sin \pi \beta}}\sum_{n=1}^\infty {{(-1)^n n \sin(\pi \beta n)} \over
{\pi^2 a^2+a^2}}.$$
The homogeneous solutions are $e^{\pm \gamma a}$ with Wronskian $W=-2\gamma$.  Then variation of
parameters, for instance, may be used to obtain the solution of the differential equation.
Otherwise, a direct verification of the stated result as satisfying the differential equation may be
performed.

\medskip
{\bf 4.123.6}
\medskip
This entry states that for Re $p>0$,
$$\int_0^\infty {{\sin ax \sinh bx} \over {\cos 2ax+\cosh 2bx}}x^{p-1}dx
={{\Gamma(p)} \over {(a^2+b^2)^{p/2}}}\sin\left(p \tan^{-1}{a \over b}\right)\sum_{k=0}^\infty {{(-1)^k}
\over {(2k+1)^p}}$$

{\it Method 1}.
We first apply the expansion of Lemma 6a, so that
$$\int_0^\infty {{\sin ax \sinh bx} \over {\cos 2ax+\cosh 2bx}}x^{p-1}dx
=-\sum_{n=1}^\infty (-1)^n\int_0^\infty e^{-2bnx}{{\sin(2nax)} \over {\cos ax}}\sinh bx ~x^{p-1}dx$$
$$=-2\sum_{n=1}^\infty \sum_{j=0}^\infty (-1)^n(-1)^j \int_0^\infty e^{-2bnx}\sin(2nax)\sinh bx 
~e^{-(2j+1)iax} x^{p-1}dx$$
$$=-2\sum_{n=1}^\infty \sum_{j=0}^\infty {{(-1)^n} \over n^p}(-1)^j \int_0^\infty e^{-2bx}\sin(2ax)
\sinh \left({b \over n}x\right) e^{-(2j+1)iax/n} x^{p-1}dx.$$
The integrand consists of exponential functions linear in $x$ together with the factor $x^{p-1}$
so that the integral evaluates in terms of $\Gamma(p)$:
$${1 \over n^p}\int_0^\infty e^{-2bx}\sin(2ax) \sinh \left({b \over n}x\right) e^{-(2j+1)iax/n} x^{p-1}dx$$
{\small
$$={{\Gamma(p)} \over 2}\left\{-{1 \over {[(2j+1)ia+(2n+1)b]^p}}\left[1+{{4a^2n^2} \over
{[(2j+1)ia+(2n+1)b]^2}}\right]^{-p/2}\sin\left[p\tan^{-1}{{2an} \over {(2j+1)ia+(2n+1)b}}\right]\right.$$
$$+\left.{1 \over {[(2j+1)ia+(2n-1)b]^p}}\left[1-{{4a^2n^2} \over
{[(2j+1)ia+(2n-1)b]^2}}\right]^{-p/2}\sin\left[p\tan^{-1}{{2an} \over {(2j+1)ia+(2n-1)b}}\right]\right
\}.$$ }
For manipulations concerning the factors $\sin(p \tan^{-1}y)$ we refer to (2.1) and (4.1) in the
following.
Performing the summation over $n$ gives the result. 

More expedient is {\it Method 2}.  We note that
$$\cosh2ax +\cosh 2bx=2\cosh(a+b)x \cosh(a-b)x,$$
and $\sinh ax\sinh bx=(1/4)[e^{(a+b)x}+e^{-(a+b)x}-e^{-(a-b)x}-e^{(a-b)x}]$, leading to
$${{\sinh ax \sinh bx} \over {\cosh 2ax +\cosh 2bx}}={1 \over 2}\left[{1 \over {e^{(a-b)x}+e^{-(a-b)x}}}
-{1 \over {e^{(a+b)x}+e^{-(a+b)x}}}\right].$$
Then
$$\int_0^\infty x^{p-1}{{\sinh ax \sinh bx} \over {\cosh 2ax +\cosh 2bx}}dx={1 \over 2}
\sum_{j=0}^\infty (-1)^j \int_0^\infty x^{p-1}\left[e^{-(2j+1)(a-b)x}-e^{-(2j+1)(a+b)x}\right]dx$$
$$={{\Gamma(p)} \over 2}\left[{1 \over {(a-b)^p}}-{1 \over {(a+b)^p}}\right]\sum_{j=0}^\infty {{(-1)^j}
\over {(2j+1)^p}}$$
$$={{\Gamma(p)} \over {2(a^2-b^2)^{p/2}}}\left[\left({{a+b} \over {a-b}}\right)^{p/2}-\left({{a+b} \over {a-b}}\right)^{-p/2}\right]\sum_{j=0}^\infty {{(-1)^j}\over {(2j+1)^p}}.$$
We may note that using (2.1)
$$ip \tan^{-1}\left({a \over b}\right)={p \over 2}\ln\left({{b+ia} \over {b-ia}}\right).$$
Then
$$\sin\left(p \tan^{-1}\left({a \over b}\right)\right)={1 \over {2i}}\left[\left({{b+ia} \over {b-ia}}\right)^{p/2}-\left({{b+ia} \over {b-ia}}\right)^{-p/2}\right]. \eqno(4.1)$$
Now let $a \to ai$, so that $\cosh2ax \to \cos 2ax$, $\sinh ax \to i \sin ax$, and
$(a^2-b^2)^{-p/2} \to (-1)^{p/2}(a^2+b^2)^{-p/2}$, and the stated result obtains.

{\it Remark}.  We may easily observe that the factor
$$\sum_{j=0}^\infty {{(-1)^j}\over {(2j+1)^p}}={1 \over 4^p}\left[\zeta\left(p,{1 \over 4}\right)-
\zeta\left(p,{3 \over 4}\right)\right].$$
Therefore, by taking derivatives with respect to $p$, and evaluating differences at $p=1$, we may 
obtain various expressions for the difference of the particular Stieltjes constants $\gamma_k(1/4)-\gamma_k(3/4)$ \cite{coffeystconstrefs}, but this is omitted.

\medskip
{\bf 4.123.7}
\medskip
This entry states that for $a>0$,
$$\int_0^\infty \sin ax^2 {{\sin(\pi x/2)\sinh(\pi x/2)} \over {\cos \pi x+\cosh \pi x}}x dx
={1 \over 4}\left[{{\partial \theta_1(z,q)} \over {\partial z}}\right]_{z=0,q=e^{-2a}}.$$

The theta function in question is given by (\cite{grad}, p.\ 921)
$$\theta_1(u,q)=2\sum_{n=1}^\infty (-1)^{n+1}q^{(n-1/2)^2}\sin(2n-1)u.$$
Therefore,
$$\left[{{\partial \theta_1(u,q)} \over {\partial u}}\right]_{u=0,q}=2\sum_{n=1}^\infty (-1)^{n+1}(2n-1)q^{(n-1/2)^2}.$$

We may use the representation
$$q^{(n-1/2)^2}=\int_0^\infty {{x\sin ax} \over {(n-1/2)^4+x^2}}dx, ~~q=\exp(-a),$$
to write the summation
$$\sum_{n=1}^\infty (-1)^{n+1}(2n-1)q^{4(n-1/2)^2}=\sum_{n=1}^\infty (-1)^{n+1}(2n-1) \int_0^\infty 
{{x\sin ax} \over {(n-1/2)^4+x^2}}dx, ~~q=\exp(-4a).$$
The summation and integration on the right side are interchanged, using the decomposition
$$\sum_{n=1}^\infty {{(-1)^{n+1}(2n-1)} \over {(n-1/2)^4+x^2}}={1 \over {2ix}}\sum_{n=1}^\infty
(-1)^{n+1}(2n-1)\left[{1 \over {(n-1/2)^2-ix}}-{1 \over {(n-1/2)^2+ix}}\right],$$
and
$$\sum_{n=1}^\infty {{(-1)^{n+1}(2n-1)} \over {(2n-1)^2+y^2}}={\pi \over {4\cosh(\pi y/2)}}.$$
This latter identity follows from either the use of the partial fractions decomposition of
the sech function, or as a special case of Fourier series.  We have the relation
$$\cosh(\pi \sqrt{i}y/2)=\cosh[\pi(1+i)y/2\sqrt{2}]=\cosh(\pi x/2\sqrt{2})\cos(\pi x/2\sqrt{2})
+i\sinh(\pi y/2\sqrt{2})\sin(\pi y/2\sqrt{2}),$$
yielding, with $y \to \pm \sqrt{x}$, the result
$$\sum_{n=1}^\infty (-1)^{n+1}(2n-1)q^{4(n-1/2)^2}={1 \over 2}\int_0^\infty {{\sinh\left({\pi \over 2}
\sqrt{x/2}\right)\sin\left({\pi \over 2}\sqrt{x/2}\right)} \over {\cosh(\pi \sqrt{x/2})+\cos(\pi\sqrt{x/2})
}}\sin a x dx.$$
Finally in the integral $x$ is replaced by $2x^2$ to give the entry.

\medskip
\centerline{\bf 5.  Concerning section 4.124}
\medskip

This section has a different `flavor' as Bessel functions of the first kind $J_n$ appear in
intermediate steps, and the final result is given in terms of the zeroth order Bessel function $J_0$.

{\bf 4.124.1}  We write this entry in a slightly extended form with upper limit $u$, not
necessarily $1$,
$$\int_0^u \cos px \cosh(q\sqrt{u^2-x^2}){{dx} \over \sqrt{u^2-x^2}}={\pi \over 2}J_0(\sqrt{p^2-q^2}u).$$
For convenience we assume $p >0$.

{\it Method 1}.
By using the Maclaurin series for the $\cosh$ function, we have
$$\int_0^u \cos px \cosh(q\sqrt{u^2-x^2}){{dx} \over \sqrt{u^2-x^2}}=\sum_{n=0}^\infty {q^{2n} \over
{(2n)!}}\int_0^u \cos px (u^2-x^2)^{n-1/2}dx$$
$$=\sqrt{\pi}\sum_{n=0}^\infty {q^{2n} \over {(2n)!}}2^{n-1}{u^n \over p^n}\Gamma(n+1/2)J_n(pu),$$
wherein the integral representation (\cite{grad}, 8.411.8, p.\ 953) has been used.
The summation is rewritten with the use of the duplication formula of the Gamma function (\cite{grad},
p.\ 938)
$$(2n)!=\Gamma(2n+1)={2^{2n} \over \sqrt{\pi}}\Gamma(n+1/2)\Gamma(n+1),$$
giving
$$\int_0^u \cos px \cosh(q\sqrt{u^2-x^2}){{dx} \over \sqrt{u^2-x^2}}=
{\pi \over 2}\sum_{n=0}^\infty {q^{2n} \over 2^n}\left({u \over p}\right)^n {{J_n(pu)} \over {n!}}.$$
Now the summation identity (\cite{grad}, 8.515.1, p.\ 974)
$$\sum_{k=0}^\infty {{(-1)^k} \over {k!}}t^k\left({{2z+t} \over {2z}}\right)^kJ_{k+\nu}(z)
=J_\nu(z+t) \eqno(5.1)$$
is used with $\nu=0$,
$$\sum_{k=0}^\infty {{(-1)^k} \over {k!}}t^k\left({{2z+t} \over {2z}}\right)^kJ_k(z)=J_0(z+t).$$
With $z=pu$, we have the quadratic equation
$$(-t)\left({{2z+t} \over {2z}}\right)={q^2 \over 2}{u \over p}$$
determining $t$, so that $t=-pu \pm \sqrt{p^2-q^2}u=-z \pm \sqrt{p^2-q^2}u$.  The choice of sign
is immaterial as $J_0(w)=J_0(-w)$, and the result follows.

{\it Remark}.  This entry may be further extended to
$$\int_0^u \cos px \cosh(q\sqrt{u^2-x^2}){{dx} \over {(u^2-x^2)^{\nu+1/2}}}
=\pi\sum_{n=0}^\infty {q^{2n} \over 2^{n+\nu+1}}{u^{n-\nu} \over p^{n-\nu}}{{\Gamma(n+1/2-\nu)} \over
{\Gamma(n+1/2)}}{{J_{n-\nu}(pu)} \over {n!}}.$$
For particular cases for $\nu$, including integer and half-integer values, the summation may be performed in terms of products or sums of Bessel functions.  
We note that (5.1) may be repeatedly differentiated with respect to $t$ to provide a family of
summation identities.  

As an example for $\nu=-1$, with $\Gamma(n+3/2)/\Gamma(n+1/2)=(n+1/2)$, using one derivative of
(5.1) with respect to $t$ and (\cite{grad}, p.\ 967), 
$(d/dz)J_{n+\nu}(z)=(1/z)[J_{n+\nu-1}(z)-(n+\nu)J_{n+\nu}(z)]$,
we obtain 
$$\int_0^u \cos px\cosh(q\sqrt{u^2-x^2})\sqrt{u^2-x^2}dx$$
$$={{\pi u^2 q^2} \over {2(q^2-p^2)}}\left[J_0(\sqrt{p^2-q^2}u)-{{(p^2+q^2)} \over {uq^2\sqrt{p^2-q^2}}}
J_1(\sqrt{p^2-q^2}u)\right].$$

The following {\it Method 2} is due to Paul Martin.  It begins by using the evenness of the
integrand, the change of variable $x=u\cos \theta$, and the addition formula for the cosine function.
$$I=\int_0^u \cos px \cosh(q\sqrt{u^2-x^2}){{dx} \over \sqrt{u^2-x^2}}
={1 \over 2}\int_{-u}^u \cos px \cosh(q\sqrt{u^2-x^2}){{dx} \over \sqrt{u^2-x^2}}$$
$$={1 \over 2}\int_0^\pi \cos(pu\cos \theta)\cosh(qu\sin \theta)d\theta
={1 \over 4}\int_{-\pi}^\pi \cos(pu\cos \theta)\cosh(qu\sin \theta)d\theta$$
$$={1 \over 8}(I_++I_-),$$
where 
$$I_{\pm}\equiv \int_\pi^\pi \cos(pu\cos\theta\pm iqu\sin\theta)d\theta.$$
By rewriting $p$ and $q$ as $p=\sqrt{p^2-q^2}\cos \alpha$ and  $iq=\sqrt{p^2-q^2}\sin \alpha$,
these integrals may be expressed as
$$I_{\pm}=\int_{-\pi}^\pi \cos [u\sqrt{p^2-q^2}\cos(\theta \mp \alpha)d\theta.$$
Now there is the expansion (\cite{nbs} 9.1.44, p.\ 361)
$$\cos(z\cos \phi)=J_0(z)+2\sum_{k=1}^\infty (-1)^k J_{2k}(z)\cos 2k\phi,$$
so that
$$\cos[z\cos (\theta \mp \alpha)]]=J_0(z)+2\sum_{k=1}^\infty (-1)^k J_{2k}(z)[\cos 2k\theta\cos2k\alpha\pm
\sin 2k\theta\sin 2k\alpha].$$
This gives
$$I_+=I_-=2\pi J_0(u\sqrt{p^2-q^2}),$$
and thus
$$I={\pi \over 2}J_0(u\sqrt{p^2-q^2}).$$

{\bf 4.124.2} This entry states
$$\int_u^\infty \cos ax\cosh\sqrt{\beta(u^2-x^2)}{{dx} \over \sqrt{u^2-x^2}}={\pi \over 2}J_0\left({u
\over \sqrt{a^2-\beta^2}}\right).$$

This entry is puzzling in several respects, as it is not clear why the integrand is not instead
written in terms of $\cos \sqrt{\beta(x^2-u^2)}$ in order to ensure convergence.  There are no conditions stated between $a$ and $\beta$,
and it is further not clear that the result should depend upon $\beta^2$ and not simply $\beta$.
In addition, just the $\beta=0$ case does not seem to agree with the stated result.
For the very special case that $a=\beta=0$, the integral appears to be divergent, while the
stated result would apparently give $0$.  
The source of this entry is cited as \cite{erdelyi} (p.\ 34), but that cosine transform is what we
have proved for the previous entry 4.124.1.
It appears that the citations for entries 4.124.1 and 4.124.2 may have been reversed in some way,
as manipulation of the exponential Fourier transform that is the first entry on p.\ 119 of 
\cite{magnus} may be the basis of an entry with an integral closely related to 4.124.2, perhaps
$$\int_a^\infty \cos b\sqrt{x^2-a^2}{{\cos xy} \over \sqrt{x^2-a^2}}dx.$$ 



\medskip
\centerline{\bf 6.  Concerning other hyperbolic-trigonometric integrals}
\medskip

Let $\tilde{\omega}$ be the Beta function value
$$\tilde{\omega}=2\int_0^1 {{dx} \over \sqrt{1-x^4}}={1 \over 2}B\left({1 \over 4},{1 \over 2}\right)
={\sqrt{\pi} \over 2}{{\Gamma(1/4)} \over {\Gamma(3/4)}}.$$
We have separately investigated the Weierstrass $\wp$ function having periods $\tilde{\omega}$ and $\tilde{\omega}i$ \cite{mwchnumbers}.  In the Laurent expansion of this function about the origin occur
the Hurwitz numbers $\tilde{H}_n$ \cite{hurwitz}.  These rational numbers for this particular
doubly-periodic function are highly analogous to the Bernoulli numbers which occur in a large 
number of expansions of singly-periodic functions.  Like the Bernoulli numbers, the Hurwitz numbers
possess many and important number theoretic properties.  While the Bernoulli numbers $B_{2k}$ are
nonzero, it is the Hurwitz numbers $\tilde{H}_{4n}$ which are nonzero.

We present a few examples of hyperbolic-trigonometric integrals that have been determined from
certain summatory properies of the Hurwitz numbers.  The last example explicitly has the period
$\tilde{\omega}$ appearing in the evaluation.  Further details and many other results are contained in
\cite{mwchnumbers}.
$${1 \over 2}\int_0^\infty {{(\cos t+1)} \over {\cosh t-\cos t}}\left[\sinh\left({t \over 2}\right)-\sin
\left({t \over 2}\right)\right]dt=1-{\pi \over 4},$$
$$\int_0^\infty {{(\cos t+1)t^2} \over {\cosh t-\cos t}}\left[\sinh\left({t \over 2}\right)+\sin
\left({t \over 2}\right)\right]dt=16,$$
and
$$\int_0^\infty {{(\cos t+1)t} \over {\cosh t-\cos t}}\left[\cosh\left({t \over 2}\right)-\cos
\left({t \over 2}\right)\right]dt=\tilde{\omega}^2-4.$$

\pagebreak
\centerline{\bf Appendix}
\medskip

Let $\Gamma$ denote the Gamma function and $\zeta$ the Riemann zeta function.

{\bf 3.527.3}  For Re $a>0$ and Re $\mu>0$,
$$\int_0^\infty {x^{\mu-1} \over {\cosh^2 ax}}dx={4 \over {(2a)^\mu}}(1-2^{2-\mu})\Gamma(\mu)\zeta(\mu-1).$$

This is equivalent to the well known representation for Re $s>0$
$$\zeta(s)={{(1-2^{1-s})^{-1}} \over {2\Gamma(s+1)}}\int_0^\infty {{t^s e^t dt} \over {(e^t+1)^2}}.$$
We have $(e^t+1)^2=4e^t \cosh^2(t/2)$ and then the change of variable $t=2ax$ is used.
Putting $s=\mu-1$ gives the result.

{\bf 3.532.1}
$$\int_0^\infty {{x^n dx} \over {a\cosh x+b \sinh x}}={2 \over {(a+b)}}\int_0^\infty {{e^{-x}x^n dx} \over
{\left[1+\left({{a-b} \over {a+b}}\right)e^{-2x}\right]}}$$
$$={2 \over {(a+b)}}\int_0^\infty\sum_{m=0}^\infty (-1)^m \left({{a-b} \over {a+b}}\right)^m x^n 
e^{-(2m+1)x}dx$$
$$={{\Gamma(2n+1)} \over {(a+b)}}\sum_{m=0}^\infty {1 \over {(2m+1)^{n+1}}}\left({{a-b} \over {a+b}}\right)^m.$$
This holds for Re a$>0$ and Re $b>0$.  As long as $n>-1$, $n$ need not be a nonnegative integer. 

Finally, we mention the large amount of redundancy within section 4.117, and between this section
and section 4.115.  Differentiation of entry 4.117.1 with respect to $a$ yields 4.117.5.
Accordingly, we have the following list of pairs of equivalent entries within this section:
(4.117.1,4.117.5), (4.117.2,4.117.6), (4.117.3,4.117.7), and (4.117.4.117.8).  The
complementary integral to entry 4.117.9 is
$$\int_0^\infty {{\sin ax} \over {1+x^2}}\coth\left({\pi \over 4}x\right)dx
=-{\pi \over 2}e^{-a}+2\cosh a \tan^{-1}(e^{-a})+\sinh a\ln \coth\left({a \over 2}\right).$$
The right side of this expression properly limits to $0$ as $a \to 0$.

Special cases of the entries of section 4.115 give several of those of 4.117, and we provide a
number of examples.

With the change of variable $x \to x/2$, putting $\beta=\pi$ and using the replacement $a\to 2a$
in 4.115.7 yields 4.117.1.  Entry 4.115.8 with $\beta=1$ and $\gamma=\pi$ gives entry 4.117.3.
Likewise, entry 4.115.9 with $\beta=\pi$ gives 4.117.3 too.  Entry 4.115.10 with $\beta=\pi/2$
gives 4.117.4, and entry 4.115.11 with $\beta=\pi$ and $b=1$ provides 4.117.7.  Again, 4.117.7
follows from 4.115.12 with $\beta=\pi$.  Entry 4.115.13 with $\beta=\pi/2$ gives 4.117.8.  

As in section 4.117, there are many equivalent pairs of entries in section 4.115.  These include
the equivalence of 4.115.9 and 4.115.11 when $b=1$ in the latter entry.

\bigskip
\centerline{\bf Acknowledgements}
Dr.\ Paul Martin is thanked for his reading of, and contributions to, section 5, including his method
of proof for entry 4.124.1.

\pagebreak

\end{document}